\documentclass[12pt,amssymb,amstex]{article}

\usepackage{amsmath}
\usepackage{amssymb}
\usepackage{amsxtra}
\usepackage{latexsym}
\usepackage{ifthen}

\usepackage[dvips]{graphicx}
\usepackage{epstopdf}
\usepackage{epsfig}
\usepackage{epic}
\usepackage{eepic}

\def\Xint#1{\mathchoice
   {\XXint\displaystyle\textstyle{#1}}
   {\XXint\textstyle\scriptstyle{#1}}
   {\XXint\scriptstyle\scriptscriptstyle{#1}}
   {\XXint\scriptscriptstyle\scriptscriptstyle{#1}}
   \!\int}
\def\XXint#1#2#3{{\setbox0=\hbox{$#1{#2#3}{\int}$}
     \vcenter{\hbox{$#2#3$}}\kern-.5\wd0}}
\def\dashint{\Xint-}

\begin{document}

\newcounter{lemma}[section]
\newcommand{\lemma}{\par \refstepcounter{lemma}%
{\bf Lemma \arabic{section}.\arabic{lemma}.}}
\renewcommand{\thelemma}{\thesection.\arabic{lemma}}

\newcounter{corol}[section]
\newcommand{\corol}{\par \refstepcounter{corol}%
{\bf Corollary \arabic{section}.\arabic{corol}.}}
\renewcommand{\thecorol}{\thesection.\arabic{corol}}

\newcounter{rem}[section]
\newcommand{\rem}{\par \refstepcounter{rem}%
{\bf Remark \arabic{section}.\arabic{rem}.}}
\renewcommand{\therem}{\thesection.\arabic{rem}}

\newcounter{theo}[section]
\newcommand{\theo}{\par \refstepcounter{theo}%
{\bf Theorem \arabic{section}.\arabic{theo}.}}
\renewcommand{\thetheo}{\thesection.\arabic{theo}}

\newcounter{propo}[section]
\newcommand{\propo}{\par \refstepcounter{propo}%
{\bf Proposition \arabic{section}.\arabic{propo}.}}
\renewcommand{\thepropo}{\thesection.\arabic{propo}}

\numberwithin{equation}{section}

\newcommand{\osc}{\operatornamewithlimits{osc}}

\def\Xint#1{\mathchoice
   {\XXint\displaystyle\textstyle{#1}}%
   {\XXint\textstyle\scriptstyle{#1}}%
   {\XXint\scriptstyle\scriptscriptstyle{#1}}%
   {\XXint\scriptscriptstyle\scriptscriptstyle{#1}}%
   \!\int}
\def\XXint#1#2#3{{\setbox0=\hbox{$#1{#2#3}{\int}$}
     \vcenter{\hbox{$#2#3$}}\kern-.5\wd0}}
\def\dashint{\Xint-}

\def\cc{\setcounter{equation}{0}
\setcounter{figure}{0}\setcounter{table}{0}}

\overfullrule=0pt

\title{{\bf About some mappings in $\bf{\lambda(r)}$-regular metric spaces}}

\author{\bf     R. Salimov,  O. Afanas'eva}

\date{\today \hskip 4mm ({\tt BOUNDARY-RSSY311211.tex})}

\maketitle

\abstract  It is formulated conditions on functions $Q(x)$ and
 boundaries of domains under which every  $Q$-homeo\-morphism admits
continuous or homeomorphic extension to the boundary in metric
spaces with measures.
\endabstract

\medskip



\large

\section{Introduction}

Mapping theory started in the 18th century. Beltrami, Caratheodory,
Christoffel, Gauss, Hilbert, Liouville, Poincar\'{e}, Riemann,
Schwarz, and so on all left their marks in this theory. Conformal
mappings and their applications to potential theory, mathematical
physics, Riemann surfaces, and technology played a key role in this
development.
\medskip

During the late 1920s and early 1930s, Gr\"otzsch, Lavrentiev, and
Morrey introduced a more general and less rigid class of mappings
that were later named quasicon\-for\-mal. The concept of
$Q$-homeo\-morphism
 is a natural extension of the geometric definition of
quasiconformality; see, e.g., \cite{MRSY}. The subject of
$Q$-homeomorphisms is interesting on its own right and has
applications to many classes of mappings. In particular, the theory
of $Q$-homeomorphisms can be applied to mappings in local Sobolev
classes (see, e.g., Sections 6.3 and 6.10 in \cite{MRSY}) to the
mappings with finite length distortion (see Sections 8.6 and 8.7 in
\cite{MRSY}) and to the finitely bi-Lipschitz mappings; see Section
10.6 in \cite{MRSY}. The main goal of the theory of
$Q$-homeomorphisms is to clear up various interconnections between
properties of the majorant $Q(x)$ and the corresponding properties
of the mappings themselves.
\medskip

\medskip

\section{Preliminaries}
Given a set $S$ in $(X,d)$ and $\alpha \in[0,\infty),$ $H^{\alpha}$
denotes the {\bf $\alpha$-dimen\-sio\-nal Hausdorff measure} of $S$
in $(X,d),$ i.e.,
\begin{equation}\label{eq333.2.1}H^{\alpha}(S)\ =\ \sup_{\varepsilon>0}\
H^{\alpha}_{\varepsilon}(S)\,,\end{equation}
\begin{equation}\label{eq333.2.2} H^{\alpha}_{\varepsilon}(S)\ =\
\inf\ \sum^{\infty}_{i=1}\delta_i^{\alpha}\,,\end{equation} where
the infimum is taken over all countable collections of numbers
$\delta_i\in(0,\varepsilon)$ such that some sets $S_i$ in $(X,d)$
with diameters $\delta_i$ cover $S.$ Note that $H^{\alpha}$ is
nonincreasing in the parameter $\alpha.$ The {\bf Hausdorff
dimension} of $S$ is the only number $\alpha\in[0,\infty]$ such that
$H^{\alpha'}(S)=0$ for all $\alpha'>\alpha$ and
$H^{\alpha''}(S)=\infty$ for all $\alpha''<\alpha .$
\medskip

Recall, for a given continuous path $\gamma: [a,b]\to X$ in a metric
space $(X,d),$ that its {\bf length} is the supremum of the sums
$$
\sum^{k}_{i=1} d(\gamma(t_i),\gamma(t_{i-1}))
$$
over all partitions $ a=t_0\leq t_1\leq\dots \leq t_k=b$ of the
interval $ [a,b].$ The path $\gamma$ is called {\bf rectifiable} if
its length is finite.
\medskip

In what follows, $(X,d,\mu)$ denotes a space $X$ with a metric $d$
and a locally finite Borel measure $\mu.$ Given a family of paths
$\Gamma$ in $X$, a Borel function $\varrho:X\to[0,\infty]$ is called
{\bf admissible} for $\Gamma$, abbr. $\varrho\in
\mathrm{adm}\,\,\Gamma$, if \begin{equation}\label{eq13.2}
\int\limits_{\gamma}\varrho\,ds\ \geq\ 1\end{equation} for all
$\gamma\in\Gamma$.
\medskip

An open set in $X$  whose points can all be connected pairwise  by
continuous paths is called a {\bf domain} in $X.$ Let $D$ and $D'$
be domains with finite Hausdorff dimensions $\alpha$ and $\alpha'\ge
1$ in spaces $(X,d,\mu)$ and $(X',d',\mu'),$ and let
$Q:D\to[0,\infty]$ be a measurable function. We say that a
homeomorphism $f:D\to{D'}$ is a {\bf $Q$-homeomorphism} if
\begin{equation}\label{eq13.4} M(f\Gamma)\leq\int\limits_{D}
Q(x)\cdot\varrho^{\alpha}(x)\ d\mu(x)\end{equation} for every family
$\Gamma$ of paths in $D$ and every admissible function $\varrho$ for
$\Gamma$.
\medskip

The {\bf modulus} of the path family $\Gamma$ in $D$ is given by the
equality \begin{equation}\label{eq13.5}M(\Gamma)\ =\
\inf\limits_{\varrho\in
\mathrm{adm}\,\Gamma}\int\limits_{D}\varrho^{\alpha}(x)\
d\mu(x).\end{equation} In the case of the path family
$\Gamma'=f\Gamma$, we take the Hausdorff dimension $\alpha'$ of the
domain $D'.$
\medskip

A space  $(X,d,\mu)$ is called {\bf (Ahlfors) $\bf{\alpha}$-regular}
if there is a constant $C\ge 1$ such that
\begin{equation}\label{eq13.6} C^{-1}r^{\alpha}\leq\mu(B_{r})\leq Cr^{\alpha} \end{equation}
for all balls $B_r$ in $X$ with the radius $r<\mathrm{diam}\,X.$ A
space $(X,d,\mu)$ is {\bf (Ahlfors) regular} if it is (Ahlfors)
$\alpha$-regular for some $\alpha\in(1,\infty)$.

A space  $(X,d,\mu)$ is {\bf upper ${\bf{\alpha}}$-regular at a
point} $x_0\in X$ if there is a constant $C> 0$ such that
\begin{equation}\label{eq3.7}
\mu(B(x_0,r))\leq Cr^{\alpha}\end{equation} for the balls $B(x_0,r)$
centered at $x_0\in X$ with all radii  $r<r_0$ for some $r_0>0.$ A
space  $(X,d,\mu)$ is {\bf upper $\bf{\alpha}$-regular} if
condition (\ref{eq3.7}) holds at every point $x_0\in X$, see
\cite{MRSY}, p. 258.

Recall that a topological space is {\bf connected space} if it is
impossible to split it into two non-empty open sets. Compact
connected spaces are called {\bf continua}. A topological space $T$
is said to be {\bf path connected} if any two points $x_1$ and $x_2$
in $T$ can be joined by a path $\gamma :[0,1]\rightarrow
T,\,\gamma(0)=x_1 $ and $\gamma(1)=x_2 $. A {\bf domain} in $T$ is
an open path connected set in $T.$  A domain $D$ in a topological
space $T$ is called {\bf locally connected at a point}
$x_0\in\partial D$ if, for every neighborhood $U$ of the point
$x_0,$ there is its neighborhood $V\subseteq U$ such that $V\cap D$
is connected, \cite{Ku2}, c. 232. Similarly, we say that a domain
$D$ is {\bf locally path connected at a point} $x_0\in\partial D$
if, for every neighborhood $U$ of the point $x_0,$ there is its
neighborhood $V\subseteq U$ such that $V\cap D$ is path connected.

The boundary of $D$ is {\bf weakly flat at a point} $x_0\in \partial
D$ if, for every number $P>0$ and every neighborhood $U$ of the
point $x_0,$ there is its neighborhood $V\subset U$ such that
\begin{equation}\label{eq13.15} M(\Delta(E,F; D))\geq P \end{equation}
for all continua $E$ and $F$ in $D$ intersecting $\partial U$ and
$\partial V.$

The boundary of the domain $D$ is {\bf strongly accessible at a
point $x_0\in
\partial D$}, if, for every neighborhood  $U$ of the point $x_0$,
there is a compact set $E\subset D$, a neighborhood $V\subset U$ of
the point $x_0$ and a number $\delta>0$ such that
\begin{equation}\label{eq13.515} M(\Delta(E,F; D))\geq \delta \end{equation}
for every continuum $F$ in $D$ intersecting $\partial U$ and
$\partial V.$

Finally, we say that the boundary $\partial D$  is {\bf weakly flat}
and {\bf strongly accessible} if the corresponding properties hold
at every point of the boundary, see \cite{MRSY}.

We start first from the following general statement, see Lemma 13.3
 and Theorem 13.3 in \cite{MRSY}.

\begin{lemma}\label{lem4} {\it Let a space $X$ be path connected at a
point $x_0\in D$ that has a compact neighborhood, let $X'$ be a
compact weakly flat space, and let $f:D\to D'$ be a
$Q$-homeomorphism, where $Q: D\to[0,\infty]$ is a measurable
function satisfying the condition
\begin{equation}\label{1} \int\limits_{D\cap A(x_0,\varepsilon,\varepsilon_{0})}
Q(x)\cdot\psi^{\alpha}_{x_0,\varepsilon}(d(x,x_0))\,
d\mu(x)\,=\,o(I^{\alpha}_{x_0}(\varepsilon))\end{equation} as
$\varepsilon\to0$, $A(x_0,\varepsilon,\varepsilon_{0})\,=\,\{x\in X:
\varepsilon<d(x,x_0)<\varepsilon_{0}\},$
 where $ \varepsilon_0<\mathrm{dist}(x_0,\partial
D)$ and $\psi_{x_0,\varepsilon}(t)$ is a family of nonnegative
(Lebesgue) measurable functions on $(0,\infty)$ such that
\begin{equation}\label{2}
0<I_{x_0}(\varepsilon)\,=\,\int\limits_{\varepsilon}^{\varepsilon_0}
\psi_{x_0,\varepsilon}(t)\,dt\,<\,\infty\,, \quad
\varepsilon\in(0,\varepsilon_0)\, .\end{equation} Then $f$ can be
extended to the point $x_{0}$ by continuity in $X'$. }\end{lemma}

\begin{theo}\label{th13.6.3} Let $D$ be locally path connected
 at all its boundary points and $\overline{D}$ compact, $D'$ with a weakly flat
boundary, and let $f:D\rightarrow D'$ be a $Q$-homeomorphism with $
Q\in L_{\mu}^{1}(D)$. Then the inverse homeomorphism
$g=f^{-1}:D'\rightarrow D$ admits a continuous extension
$\overline{g}:\overline{D'}\rightarrow\overline{D}$. \end{theo}

\section{Basic results}

We will say that a space  $(X,d,\mu)$ is {\bf upper ${\bf
\lambda(r)}$-regular at a point} $x_0\in X$ if there is a constant
$C> 0$ such that
\begin{equation}\label{eq13.7}
\mu(B(x_0,r))\leq C\lambda_{x_0}(r)\end{equation} for the balls
$B(x_0,r)$ centered at $x_0\in X$ with all radii  $r<r_0$ for some
$r_0>0,$ and $\lambda(r)$ is increasing function. We will also say
that a space $(X,d,\mu)$ is {\bf upper ${\bf \lambda(r)}$-regular}
if condition (\ref{eq13.7}) holds at every point  $x_0\in X$.

\begin{lemma}\label{lem2} {\it Let $D$ be a domain in a space $(X,d,\mu)$
that is $upper \lambda(r)$-regular at the point $x_0\in
\overline{D}$ and $\lambda$ is increasing function. If for every
nonnegative function $\varphi:D\to\mathbb{R}$  condition
\begin{equation}\limsup\limits_{\varepsilon\rightarrow 0}\
\frac{1}{\mu(B(x_0,\varepsilon)\cap
D)}\int\limits_{B(x_0,\varepsilon)\cap D}
|\varphi(x)|d\mu(x)<\infty\end{equation} holds then
\begin{equation}\label{eq13.4.5}\int\limits_{D\cap A(x_0,\varepsilon,\varepsilon_{0})}
\frac{|\varphi(x)|\,d\mu(x)}{\lambda(d(x,x_0))}
=O\left(\log\frac{1}{\varepsilon}\right)\end{equation} as
$\varepsilon\to0$ and some $\varepsilon_{0}\in(0,\delta_0),$ where
$\delta_0=\min\,(e^{-1},d_0)$, $d_0=\sup\limits_{x\in D}d(x,x_0)$,
$A(x_0,\varepsilon,\varepsilon_{0})\,=\,\{x\in X:
\varepsilon<d(x,x_0)<\varepsilon_{0}\}, B_r=B(x_0,r)=\{ x\in X :
d(x,x_0)<r\}.$}\end{lemma}

{\it Proof.} Choose $\varepsilon_0\in (0,\delta_0)$ such that the
function $\varphi$ is integrable in $D$ with respect to  the measure
$\mu $, where
$$\delta_0=\sup_{r\in(0,\varepsilon_0)}\ \ \frac{1}{\mu(D_r)}\int\limits_{D_r}
|\varphi(x)|\ d\mu(x)\ <\infty,$$ $D_r=D\cap B_r.$ Further, let
$\varepsilon<2^{-1}\varepsilon_0,\
\varepsilon_k=2^{-k}\varepsilon_0,\ A_k=\{ x\in X:
\varepsilon_{k+1}\leq d(x,x_0)<\varepsilon_k\},$
$B_k=B(\varepsilon_k).$ Choose a natural number $N$ such that
$\varepsilon\in[\varepsilon_{N+1},\varepsilon_N)$. Then $D\cap
A(x_0,\varepsilon,\varepsilon_0) \subset \Delta
(\varepsilon):=\bigcup_{k=0}^{N}\Delta_{k}$, where $\Delta_{k}=D\cap
A_k$ and
$$\eta(\varepsilon)=
\int\limits_{\Delta(\varepsilon)}\frac{|\varphi(x)|d\mu(x)}{\lambda_{x_0}(d(x,x_0))}\leq\
\sum_{k=1}^N\int\limits_{\Delta_{k}}\frac{|\varphi(x)|d\mu(x)}{\lambda_{x_0}(d(x,x_0))}\leq\\
$$
$$\leq\sum\limits_{k=1}^N \frac{1}{\lambda(\varepsilon_{k+1})}
\int\limits_{B(\varepsilon_k)\cap D}\ |\varphi(x)|\
d\mu(x)\leq\delta_0 \cdot \sum_{k=1}^N\frac{\mu(B_k\cap
D)}{\lambda_{x_0}(\varepsilon_{k+1})}.
$$ In pursuance of upper  $\lambda(r)$-regularity we obtain that
$\mu(B_k)\leq C\cdot \lambda_{x_0}(\varepsilon_k)$ and
$$\eta(\varepsilon)\leq C\cdot \delta_0\cdot
\sum_{k=1}^N\frac{\lambda_{x_0}(\varepsilon_k)}{\lambda_{x_0}(\varepsilon_{k+1})}\leq
C\cdot \delta_0\cdot N.$$ Since $N<
\log_2{\frac{1}{\varepsilon}}=\frac{\log\frac{1}{\varepsilon}}{\log
2},$ see \cite{MRSY} p.266, then

$$\int\limits_{D\cap A(x_0,\varepsilon,\varepsilon_0)}
\frac{|\varphi(x)|d\mu(x)}{\lambda_{x_0}(d(x,x_0))} \leq
\frac{C\cdot \delta_0}{\log 2}\log{\frac{1}{\varepsilon}}.$$ So, we
complete the proof.

As before, here $(X,d,\mu)$ and $(X',d',\mu')$ are spaces with
metrics $d$ and $d'$ and locally finite Borel measures $\mu$ and
$\mu'\,,$ and  $D$ and $D'$ are domains in $X$ and $X'$ with finite
Hausdorff dimensions $\alpha$ and $\alpha'> 1$, respectively.

\begin{theo}\label{theo3} {\it Let $X$ be upper $\lambda(r)$-regular
at a point $x_0\in \partial D$ where $D$ is locally path connected,
$\overline{D'}$ be compact and $\partial D'$ strongly accessible
with \begin{equation}\label{56} \lambda(\varepsilon)=
o\left(\varepsilon^{\alpha}\log^{\alpha-1}\frac{1}{\varepsilon}\right)
\end{equation} as $\varepsilon\rightarrow 0.$
 If
\begin{equation}\label{6}
\limsup\limits_{\varepsilon\rightarrow
0}\frac{1}{\mu(B(x_0,\varepsilon)\cap
D)}\int\limits_{B(x_0,\varepsilon)\cap D} Q(x)
d\mu(x)<\infty,\end{equation}
 then any
$Q$-homeomorphism $f:D\to D'$ can be extended to the point $x_{0}$
by continuity in $(X',d')$.}\end{theo}

{\it Proof.} By elementary argument, we see that condition
(\ref{56}) implies that
\begin{equation}\label{3} \int\limits_{0}^{\varepsilon_0}
\frac{dt}{\lambda^{\frac{1}{\alpha}}(t)}=\infty.\end{equation}
Choosing in Lemma \ref{lem4}
$\psi(t)=\frac{1}{\lambda^{1/\alpha}(t)}$
 and combining conclusion Lemma \ref{lem2}, we show, by L'Hospital rule, that
 $$\lim\limits_{\varepsilon\rightarrow 0}
\int\limits_{D\cap A(x_0,\varepsilon,\varepsilon_{0})}
 \frac{Q(x)d\mu(x)}{\lambda(d(x,x_0))}\ \big/
\left(\int\limits_{\varepsilon}^{\varepsilon_0}
\frac{dt}{\lambda^{1/\alpha}(t)}\right)^{\alpha}\leq\lim\limits_{\varepsilon\rightarrow
0} c\log\frac{1}{\varepsilon}\
\big/{\left(\int\limits_{\varepsilon}^{\varepsilon_0}
\frac{dt}{\lambda^{1/\alpha}(t)}\right)^{\alpha}}=$$
$$=\gamma\ \lim\limits_{\varepsilon\rightarrow
0}\frac{\lambda(\varepsilon)}{\varepsilon^{\alpha}}\log^{1-\alpha}\frac{1}{\varepsilon},
$$ where $\gamma=\left(\frac{c}{\alpha}\right)^{\alpha}$. According to (\ref{56})
$$\lim\limits_{\varepsilon\rightarrow
0}\int\limits_{D\cap A(x_0,\varepsilon,\varepsilon_{0})}
 \frac{Q(x)d\mu(x)}{\lambda(d(x,x_0))}\  \big/
\left(\int\limits_{\varepsilon}^{\varepsilon_0}
\frac{dt}{\lambda^{1/\alpha}(t)}\right)^{\alpha}=0.$$ Since
conditions of Lemma \ref{lem4} are hold then exist extending to
$x_0$ by continuity.

Combining the theorems (\ref{th13.6.3}) and (\ref{theo3}), we obtain
the following theorem.

\begin{theo}\label{theo4} {\it Let $X$ be upper $\lambda(r)$-regular
at a point $x_0\in \partial D$ where $D$ and $D'$ have weakly flat
boundaries, let  $\overline{D}$ and $\overline{D'}$ be compact, and
let $Q:D\rightarrow[0,\,\infty]$ be a function of the class
$L_{\mu}^{1}(D)$ with condition (\ref{6}). If
 \begin{equation}\label{5}
\lambda(\varepsilon)=
o\left(\varepsilon^{\alpha}\log^{\alpha-1}\frac{1}{\varepsilon}\right)
\end{equation} as $\varepsilon\rightarrow 0$  then any  $Q$-homeomorphism $f:D\to D'$ is
extended to a homeomorphism $\overline{f}:\overline{D}\to
\overline{D'}$.}\end{theo}

\medskip

\noindent  Ruslan Salimov, Olena Afanas'eva\\
Institute of Applied Mathematics and Mechanics,\\
National Academy of Sciences of Ukraine, \\
74 Roze Luxemburg str., 83114 Donetsk, Ukraine \\
es.afanasjeva@yandex.ru, salimov07@rambler.ru\\

\end{document}